\documentclass{gtart}

\def\ifplaintex{\expandafter\ifx\csname documentclass\endcsname\relax}

\expandafter\ifx\csname epsfbox\endcsname\relax\input epsf\fi

\ifplaintex 
\else
\fi

\def\gt{{\mathsurround=0pt\it $\cal G\mskip-2mu$eometry \&\ 
$\cal T\!\!$opology}}        

\def\gtp{{\mathsurround=0pt\it $\cal G\mskip-2mu$eometry \&\ 
$\cal T\!\!$opology $\cal P\!$ublications}}  

\def\lognumber#1{\def\thelognumber{#1}}
\def\volumenumber#1{\def\thevolumenumber{#1}}
\def\papernumber#1{\def\thepapernumber{#1}}
\def\volumeyear#1{\def\thevolumeyear{#1}}

\def\pagenumbers#1#2{\def\startpage{#1}\def\finishpage{#2}}
\def\published#1{\def\publishdate{#1}}
\def\proposed#1{\def\theproposer{#1}}
\def\seconded#1{\def\theseconders{#1}}
\def\received#1{\def\receiveddate{#1}}
\def\revised#1{\def\reviseddate{#1}}
\def\accepted#1{\def\accepteddate{#1}}

\def\coverauthors#1{\def\thecoverauthors{#1}}
\def\asciiauthors#1{\def\theasciiauthors{#1}}
\def\asciiaddress#1{\def\theasciiaddress{#1}}

\let\\\par\let\thevolumenumber\relax\let\thepapernumber\relax
\let\thevolumeyear\relax\let\thesamplenumber\relax\let\startpage\relax
\let\finishpage\relax\let\publishdate\relax\let\receiveddate\relax
\let\reviseddate\relax\let\accepteddate\relax\let\theasciititle\relax
\let\theasciiauthors\relax\let\theasciiaddress\relax
\let\theasciiabstract\relax
\let\theasciiemail\relax\let\theshortauthors\relax\let\theshorttitle\relax
\let\thecoverauthors\relax

\long\def\maketitlep{   

\count0=\startpage

\gt\hfill      
\hbox to 77pt{\vbox to 0pt{\vglue -15pt\epsfbox{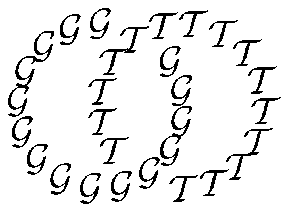}\vss}\hss}
\break
{\small\ifx\thesamplenumber\relax 
Volume \else Sample
\fi\thevolumenumber\ (\thevolumeyear)
\startpage--\finishpage\nl
Published: \publishdate}
\vglue 0.5truein plus 0.4fil minus 0.1truein

{\parskip=0pt\leftskip 0pt plus 1fil\def\\{\par\smallskip}{\ifplaintex\large
\else\Large\fi\bf\thetitle}\par\medskip}   

\vglue 0pt plus 0.1fil 

{\parskip=0pt\leftskip 0pt plus 1fil\def\\{\par}{\sc\theauthors}
\par\medskip}

\vglue 0pt plus 0.1fil 

{\small\parskip=0pt\let\newline\\
{\leftskip 0pt plus 1fil\def\\{\par}{\sl\theaddress}\par}
\expandafter\ifx\theemail\relax    
\relax\else\vglue 5pt plus 0.02fil minus 2pt\def\\{\stdspace{\rm 
and}\stdspace} 
\cl{Email:\stdspace\tt\theemail}\fi
\ifx\theurl\relax                  
\relax\else\vglue 5pt plus 0.02fil minus 2pt\def\\{\stdspace{\rm 
and}\stdspace}
\cl{URL:\stdspace\tt\theurl}\fi\par}

\vglue 7pt plus 0.3fil minus 3pt

{\bf Abstract}
\vglue 5pt plus 0.1fil minus 2pt

\theabstract

\vglue 7pt plus 0.3fil minus 3pt

{\bf AMS Classification numbers}\quad Primary:\quad \theprimaryclass

Secondary:\quad \thesecondaryclass

\vglue 5pt plus 0.3fil minus 2pt

{\bf Keywords:}\quad \thekeywords

\vglue 10pt plus 0.5fil minus 5pt

{\small  Proposed: \theproposer\hfill Received: \receiveddate\nl
Seconded: \theseconders\hfill 
\ifx\reviseddate\relax                         
Accepted: \accepteddate                        
\else
Revised: \reviseddate                          
\fi}
\eject
}       

\let\maketitlepage\maketitlep
\let\maketitle\maketitlepage

\font\phead=cmsl9 scaled 950
\font=cmsl9 scaled 1050
\font\pnum=cmbx10 scaled 913
\font\lnum=cmbx10 
\font\pfoot=cmsl9 scaled 950
\font=cmsl9 scaled 1050
\ifplaintex
\headline{\vbox to 0pt{\vskip -4.5mm\line{\small\phead\ifnum
\count0=\startpage ISSN 1364-0380 (on line)
1465-3060 (printed) \hfill {\pnum\folio}\else\ifodd\count0\def\\{ }%
\ifx\theshorttitle\relax\thetitle\else\theshorttitle\fi\hfill{\pnum\folio}
\else\def\\{ and }{\pnum\folio}\hfill\ifx\theshortauthors\relax\theauthors
\else\theshortauthors\fi\fi\fi}\vss}}
\footline{\vbox to 0pt{\vglue 0mm\line{\small\pfoot\ifnum\count0=\startpage
\copyright\ \gtp\hfill\else
\gt, Volume \thevolumenumber\ (\thevolumeyear)\hfill\fi}\vss
}}
\else
\makeatletter
\def\@oddhead{{\small\lhead\ifnum\count0=\startpage ISSN 1364-0380 (on line)
1465-3060 (printed) \hfill {\lnum\number\count0}\else\ifodd\count0
\def\\{ }\ifx\theshorttitle\relax \thetitle \else\theshorttitle\fi\hfill
{\lnum\number\count0}\else\def\\{ and }{\lnum\number\count0}
\hfill\ifx\theshortauthors\relax 
\theauthors\else\theshortauthors\fi\fi\fi}}\def\@evenhead{@oddhead}
\def\@oddfoot{\small\lfoot\ifnum\count0=\startpage\copyright\ \gtp\hfill\else
\gt, Volume \thevolumenumber\ (\thevolumeyear)\hfill\fi}
\def\@evenfoot{@oddfoot}
\makeatother
\fi

\newwrite\gtoutfile
\long\gdef\makeheadfile{  
{\def\\{, }\def\s{ }
\immediate\openout\gtoutfile head.xxx
\immediate\write\gtoutfile{To: math@arxiv.org}
\immediate\write\gtoutfile{Subject: put OR rep NNNNN:pppp}
\immediate\write\gtoutfile{--text follows this line--}
\immediate\write\gtoutfile{Proxy-for: \ifx\theasciiauthors\relax
\theauthors\else\theasciiauthors\fi\s<\ifx\theasciiemail\relax\theemail\else\theasciiemail\fi>}
\immediate\write\gtoutfile{\noexpand\\}
\immediate\write\gtoutfile{Authors: \ifx\theasciiauthors\relax
\theauthors\else\theasciiauthors\fi}
{\def\\{ }\immediate\write\gtoutfile{Title: \ifx\theasciititle\relax
\thetitle\else\theasciititle\fi}}
\immediate\write\gtoutfile{Subj-class: GT or GR or SG or ...}
\immediate\write\gtoutfile{MSC-class: \theprimaryclass\ifx\thesecondaryclass\relax\else, \thesecondaryclass\fi}
\immediate\write\gtoutfile{Journal-ref: Geom. Topol. \thevolumenumber\s
(\thevolumeyear) \startpage-\finishpage}
\immediate\write\gtoutfile{Comments: Published in Geometry and Topology at}
\immediate\write\gtoutfile{    http://www.maths.warwick.ac.uk/gt/GTVol\thevolumenumber/paper\thepapernumber.abs.html}
\immediate\write\gtoutfile{\noexpand\\}
\immediate\write\gtoutfile{}
\ifx\theasciiabstract\relax
\immediate\write\gtoutfile{\theabstract}\else
\immediate\write\gtoutfile{\theasciiabstract}\fi
\immediate\write\gtoutfile{}
\immediate\write\gtoutfile{\noexpand\\}
\immediate\write\gtoutfile{}
\immediate\closeout\gtoutfile}}  

\def\maketitlepage{\maketitlep\makeheadfile}
\let\maketitle\maketitlepage

\lognumber{260}

\volumenumber{7}\papernumber{27}\volumeyear{2003}
\pagenumbers{933}{963}
\received{5 June 2002}
\revised{4 November 2003}
\published{11 December 2003}
\accepted{5 December 2003}
\proposed{Benson Farb}
\seconded{Jean-Pierre Otal, Walter Neumann}

\usepackage[all]{xy}
\usepackage{amssymb}

\theoremstyle{definition}
\newtheorem{defi}{Definition}[section]
\newtheorem*{rem}{Remark}

\theoremstyle{plain}
\newtheorem{theo}[defi]{Theorem}
\newtheorem{lemma}[defi]{Lemma}
\newtheorem{prop}[defi]{Proposition}
\newtheorem{coro}[defi]{Corollary}

\newcommand{\dist}{{\rm dist}}
\newcommand{\G}{\Gamma}

\title{Combination of convergence groups}

\author{Fran\c{c}ois Dahmani}
\asciiauthors{Francois Dahmani}
\coverauthors{Fran\noexpand\c{c}ois Dahmani}
\address {Forschungsinstitut f\"ur Mathematik\\ETH Zentrum, R\"amistrasse, 
101\\8092 Z\"urich, Switzerland.} 
\asciiaddress {Forschungsinstitut fur Mathematik\\ETH Zentrum, Ramistrasse, 
101\\8092 Zurich, Switzerland.} 
\email{dahmani@math.ethz.ch}

\begin{abstract}
 We state and prove a combination theorem for
relatively hyperbolic groups seen as geometrically finite convergence groups. For that, we explain how to contruct a boundary for a group that is an acylindrical amalgamation of relatively hyperbolic groups over a fully quasi-convex subgroup. We apply our result to Sela's theory on limit groups and
prove their relative hyperbolicity. We also get a proof of the Howson property for limit groups.
\end{abstract}



\begin{document}

\maketitlepage

The aim of this paper is to explain how to amalgamate geometrically
finite convergence groups, or in another formulation, relatively
hyperbolic groups, and to deduce the relative hyperbolicity of Sela's
limit groups.

A group acts as a convergence group on a compact space $M$ if it acts
properly discontinuously on the space of distinct triples of $M$ (see the works
of F~Gehring,  G~Martin, A~Beardon, B~Maskit, B~Bowditch, and P~Tukia
\cite{MG}, \cite{BeM}, \cite{Bconf}, \cite{Tuk}). The convergence action is
uniform if $M$
consists only of conical limit points; the action is geometrically finite
(see \cite{BeM}, \cite{Bo1}) if $M$
consists only of conical limit points and of bounded parabolic
points. The definition of conical limit points is
a dynamical formulation of the so called points of approximation, in the
language of Kleinian groups. A point of $M$ is "bounded parabolic" is its
stabilizer acts properly discontinuously and cocompactly on its complement in
$M$, as it is the case for parabolic points of geometrically finite Kleinian
groups acting on their limit sets (see \cite{BeM}, \cite{Bo1}). See Definitions
1.1--1.3 below.

Let $\Gamma$ be a group acting properly discontinuously by isometries on a proper Gromov-hyperbolic space
$\Sigma$. Then $\Gamma$ naturally acts by homeomorphisms on the boundary
$\partial \Sigma$. If it is a uniform convergence action, $\Gamma$ is
hyperbolic in the sense of Gromov, and if the action is geometrically finite, 
following B~Bowditch \cite{Brel} we say that $\Gamma$ is
hyperbolic relative to the family $\mathcal{G}$ of the maximal
parabolic subgroups, provided that these subgroups are finitely generated. 
In such a  case, the pair $(\Gamma,
\mathcal{G})$ constitutes a relatively hyperbolic group in the sense of Gromov and Bowditch. 
Moreover, in \cite{Brel},
Bowditch explains that the compact space $\partial \Sigma$ is canonically
associated to $(\Gamma, \mathcal{G})$: it does
not depend on the choice of the space $\Sigma$. For this reason, we call it the
Bowditch boundary of the relatively hyperbolic group.

 The definitions of relative hyperbolicity in
\cite{Brel} (including the one mentionned above) are equivalent to 
Farb's relative hyperbolicity \emph{with}
the property BCP, defined in \cite{Farb} (see \cite{Szcz}, \cite{Brel}, and the
appendix of \cite{DahThesis}).

 Another theorem
of Bowditch \cite{Btop} states that the uniform convergence groups on
perfect compact spaces are exactly the hyperbolic groups acting on their Gromov
boundaries. A~Yaman \cite{Y} proved the relative version of this theorem:
geometrically finite convergence groups on perfect compact spaces with finitely generated maximal 
parabolic subgroups are exactly the
relatively hyperbolic groups acting on their Bowditch boundaries (stated below
as Theorem \ref{theo;asli}).

 We are going to formulate a definition of
quasi-convexity (Definition \ref{def;full-q-c}), generalizing an idea of
Bowditch described in \cite{Bconf}. A subgroup $H$ of a geometrically finite
convergence group on a compact space $M$ is \emph{fully quasi-convex}  if it is
geometrically finite on its limit set $\Lambda H \subset M$, and if only
finitely many translates of $\Lambda H$ can intersect non trivially together.
 We also use the notion of \emph{acylindrical} amalgamation, formulated by Sela  \cite{SelaAcyl}, which means
 that there is a number $k$ such that the stabilizer of any segment of
length $k$ in the Serre tree, is finite.

\def\Gap{$\phantom{99}$}

\begin{theo}[Combination theorem]\label{theo;comb}\Gap

{\rm(1)}\qua Let $\Gamma$ be the fundamental group of an acylindrical
finite graph of relatively hyperbolic groups, whose edge groups are
fully quasi-convex subgroups of the adjacent vertices groups. Let
$\mathcal{G}$ be the family of the images of the maximal parabolic
subgroups of the vertices groups, and their conjugates in $\Gamma$.
Then, $(\Gamma,\mathcal{G})$ is a relatively hyperbolic group.

 {\rm(2)}\qua Let $G$ be a group which is hyperbolic relative to a family of
subgroups $\mathcal{G}$, and let $ P$ be a group in  $\mathcal{G}$. Let $A$ be a
finitely generated group in which $P$ embeds as a subgroup. Then, $\Gamma = A*_P
G$ is hyperbolic relative to the family $(\mathcal{H} \cup
\mathcal{A})$, where   $\mathcal{H}$  is the set of the conjugates of the
images of elements of $\mathcal{G}$ not conjugated to $P$ in $G$, and
where $\mathcal{A}$ is the set of the conjugates of $A$ in $\Gamma$.

 {\rm(3)}\qua Let $G_1$ and $G_2$ be relatively hyperbolic groups, and let $P$ be a
maximal parabolic subgroup of $G_1$, which is isomorphic to a parabolic
(not necessarly maximal) subgroup of $G_2$. Let $\Gamma = G_1 *_P G_2$. Then
$\Gamma$ is hyperbolic relative to the family of the conjugates of the
maximal parabolic subgroups of $G_1$, except $P$, and of the conjugates of the
maximal parabolic subgroups of $G_2$.

{\rm$(3')$}\qua Let $G$ be a relatively hyperbolic group and let $P$ be a
maximal parabolic subgroup of $G$ isomorphic to a subgroup of another parabolic
subgroup $P'$ not  conjugated to $P$. Let $\Gamma = G *_P$
according to the two images. Then
$\Gamma$ is hyperbolic relative to the family of the conjugates of the
maximal parabolic subgroups of $G$, except $P$ (but including the parabolic
group $P'$).

\end{theo}

Up to our knowledge, the assumption of finite 
generation of the maximal parabolic subgroups is useful for a proof of the equivalence of different definitions 
of relative hyperbolicity. 
For the present work, it is not essential, and without major change, one can state a 
combination theorem for 
groups acting as geometrically finite convergence groups on metrisable compact spaces in general.

A first example of application of the main theorem is already known as a 
consequence of Bestvina and Feighn Combination Theorem \cite{BF}, \cite{BF2}, where there are no parabolic group: 
acylindrical amalgamations of hyperbolic groups over
quasi-convex subgroups satisfy the first case of the theorem 
(see Proposition \ref{prop;hyp}).
Another important example
is the amalgamation of relatively hyperbolic groups over a parabolic
subgroup, which is stated as the third and fourth case. They are in fact
consequences of the two first cases.

Instead of choosing the point of view of Bestvina and Feighn \cite{BF}, \cite{BF2}, 
and constructing a hyperbolic space on
which the group acts in an adequate way (see also the works of R~Gitik,
O~Kharlampovich, A~Myasnikov, and I~Kapovich, \cite{Gitik}, \cite{Kharlam},
\cite{Kapo}), we adopt a dynamical point of view: from the actions of the
vertex groups on their Bowditch's boundaries, we construct a metrizable compact
space on which $\Gamma$ acts naturally, and we check (in section 3) that this
action is of convergence and geometrically finite. At the end of the third part,
we prove the Theorem 0.1 using Bowditch--Yaman's Theorem \ref{theo;asli}.

 In other words, we construct directly the boundary of the
group $\Gamma$. This is done by gluing together the boundaries of the
stabilizers of vertices in the Bass--Serre tree, along the limit sets of the
stabilizers of the edges. This does not give a compact space, but the boundary of
the Bass--Serre tree itself naturally compactifies it. This construction is
explained in detail in section 2.

  Thus, we have a good description of the
boundary of the amalgamation. In particular:

\begin{theo}[Dimension of the boundary]\Gap

 Under the hypothesis of Theorem 0.1, let $\partial \Gamma$ be the boundary of
the relatively hyperbolic group $\Gamma$.  If the
topological  dimensions of the boundaries of the vertex groups (resp. of the
edge groups) are smaller than $r$ (resp. than $s$), then
$\hbox{dim}(\partial\Gamma) \leq \hbox{Max}\{r,s+1\}$.
\end{theo}

The application we have in mind is the study of Sela's limit groups, or equivalently $\omega$--residually
free groups \cite{SelaLimit}, \cite{KM}. In part 4, we answer the first question of Sela's list of problems
\cite{SelaProb}.

\begin{theo}
  Limit groups are hyperbolic relative to their maximal abelian non-cyclic
subgroups. 
 \end{theo}

This allows us to get some corollaries.

\begin{coro}

 Every limit group satisfies the Howson property: the intersection of two
finitely generated subgroups of a limit group is finitely generated.
\end{coro}

\begin{coro}
 Every limit group admits a $\mathcal{Z}$--structure in the sense of Bestvina
(\cite{Be}, \cite{Dah}).

\end{coro}

The first corollary was previously proved by I~Kapovich in \cite{Kapofaux}, for
\emph{hyperbolic} limit groups (see also \cite{Kapoerratum}).

  I am grateful to T~Delzant, for his interest, and advices, and to Z~Sela who
suggested the problem about   limit groups to me. I also want to thank
B~Bowditch, I~Kapovich, G~Swarup, and F~Paulin for their comments and
questions. Finally I am deeply grateful to the referee for his/her remarks.

\section{Geometrically finite convergence groups, and relative hyperbolicity}

\subsection{Definitions}

We recall the definitions of \cite{BeM}, \cite{Bconf} and \cite{Tuk}.

\begin{defi}(Convergence groups)\label{def;convergence}

A group $\Gamma$ acting on a metrizable compact space $M$ is a \emph{convergence
group} on $M$ if it acts properly discontinuously on the space of distinct
triples of $M$.
\end{defi}

If the compact space $M$ has more than two points, this is
equivalent to say that the action is of convergence if, for any sequence $(\gamma_n)_{n\in \mathbb{N}}$ of
elements of $\Gamma$, there exists
two points $\xi$ and $\zeta$ in $M$, and a subsequence
$(\gamma_{\phi(n)})_{n\in \mathbb{N}}$, such that for any
compact subspace $K \subset M\setminus \{\xi\}$, the sequence $(\gamma_{\phi(n)}
K)_{n\in \mathbb{N}}$, uniformly converges to $\zeta$.

\begin{defi}(Conical limit point, bounded parabolic point)\label{def;conlim}

 Let $\Gamma$ be a convergence group on a metrizable compact space $M$. A point
$\xi \in M$ is a  \emph{conical limit point} if there exists a sequence in
$\Gamma$, $(\gamma_n)_{n\in \mathbb{N}}$,  and two points $\zeta \neq \eta$, in
$M$, such that $\gamma_n \xi \to \zeta$ and $\gamma_n \xi' \to \eta$  for all
$\xi' \neq \xi$.

 A subgroup $G$ of $\Gamma$ is parabolic if it is infinite, fixes a point
$\xi$, and contains no loxodromic element (a loxodromic element is an element of infinite order fixing exactly two points in the boundary). In this case, the fixed point of $G$ is unique and is referred to as a parabolic point. Such a point $\xi \in M$ is
\emph{bounded parabolic} if its stabilizer $Stab(\xi)$ acts properly discontinuously
co-compactly  on $M\setminus \{\xi\}$. \end{defi}

Note that the stabilizer of a parabolic point is a maximal
parabolic subgroup of $\Gamma$.

\begin{defi}(Geometrically finite groups)\label{def;gf}

A convergence group on a compact space $M$ is \emph{geometrically finite}
if $M$ consists only of conical limit points and
bounded parabolic points.
\end{defi}

Here is a geometrical counterpart (see \cite{G}, \cite{Brel}).

\begin{defi}(Relatively hyperbolic groups)\label{def;RHG}

  We say that a group $\Gamma$ is \emph{hyperbolic relative to} a family of
finitely  generated subgroups $\mathcal{G}$, if it acts properly
discontinuously by isometries, on a proper hyperbolic space $\Sigma$, such
that the induced action on $\partial \Sigma$ is of convergence,  geometrically
finite, and such that the maximal parabolic subgroups are exactly the elements
of $\mathcal{G}$.

In this situation we also say that the pair $(\Gamma, \mathcal{G})$ is a
relatively hyperbolic group.

 \end{defi}

The boundary of $\Sigma$ is canonical in this case (see \cite{Brel}); we call
it the boundary of the relatively hyperbolic group $(\Gamma, \mathcal{G})$, or the
Bowditch boundary, and we write it $\partial \Gamma$.

As recalled in the introduction, one has:

\begin{theo}[Yaman \cite{Y}, Bowditch \cite{Btop} for groups without parabolic
subgroups]\label{theo;asli}\Gap

  Let $\Gamma$ be a geometrically finite convergence group on a perfect
metrizable compact space $M$, and let $\mathcal{G}$ be the family of its maximal
parabolic subgroups. Assume that each element of $\mathcal{G}$ is finitely
generated. Assume that there are only finitely many orbits of bounded parabolic
 points.  Then $(\Gamma, \mathcal{G})$ is relatively hyperbolic, and
$M$ is equivariantly homeomorphic to $\partial \Gamma$.

\end{theo}

In fact, by a result of Tukia (\cite{TukiaCrelle}, Theorem 1B),  the
assumption of finiteness of the set of orbits of parabolic points can be
omitted.  With this dictionary between geometrically finite convergence groups,
and relatively hyperbolic groups, we will sometimes say that a group $\Gamma$ is
relatively hyperbolic with Bowditch boundary $\partial \Gamma$, when we mean
that the pair $(\Gamma, \mathcal{G})$ is relatively hyperbolic, where
$\mathcal{G}$ is the family of maximal parabolic subgroups in the action on
$\partial \Gamma$.

\subsection{Fully quasi-convex subgroups}

Let $\Gamma$ be a convergence group on $M$. According to \cite{Bconf}, the \emph{limit set} $\Lambda H$ of
an infinite non virtually cyclic subgroup $H$,  is the unique minimal non-empty
closed $H$--invariant subset of $M$. The limit set of a virtually cyclic
subgroup of $\Gamma$ is the set of its fixed points in $M$, and the limit set
of a finite group is empty. We will use this for relatively hyperbolic groups
acting on their Bowditch boundaries.

\begin{defi}(Quasi-convex and fully quasi-convex subgroups)\label{def;full-q-c}

 Let $\Gamma$ be a relatively hyperbolic group, with Bowditch boundary $\partial
\Gamma$, and let $H$ be a group acting as a geometrically finite convergence group on a compact space 
$\partial H$. We assume that $H$  embeds in $\Gamma$ as a subgroup.  We
say that $H$ is \emph{quasi-convex} in $\Gamma$ if its limit set $\Lambda H
\subset \partial \Gamma$  is equivariantly homeomorphic to $\partial H$.

It is \emph{fully quasi-convex} if it is quasi-convex and if, for any infinite
sequence $(\gamma_n)_{n\in \mathbb{N}}$ all in distinct left cosets of $H$, the
intersection $\bigcap_{n} (\gamma_n \Lambda H)$ is empty.

 \end{defi}

\textbf{Remark (i)}\qua If $H$ is a subgroup of  $ \Gamma$, and if $\Gamma$ acts
as a convergence group on a compact space $M$, every conical limit point for $H$
acting on  $\Lambda H \subset M$, is a conical limit point for $H$ acting in
$M$, and therefore, even for $\Gamma$ acting on $M$. Therefore it is not a
parabolic point (see the result of Tukia, described in \cite{Bconf} Prop.3.2,
see also \cite{TukiaCrelle}), and each parabolic point for $H$ in $\Lambda H$ is
a parabolic point for $\Gamma$ in $M$, and its maximal parabolic subgroup in $H$
is exactly the intersection of its maximal parabolic subgroup in $\Gamma$ with
$H$.

\medskip

\textbf{Remark (ii)}\qua if $H$ is a quasiconvex subgroup of a relatively
hyperbolic group $\Gamma$, and if its maximal parabolic
subgroups are finitely generated, then it is hyperbolic relative 
to these maximal parabolic
subgroups (by Theorem \ref{theo;asli}), hence 
it is finitely generated. In particular, it
is always the case when the parabolic subgroups of $\Gamma$ are finitely
generated abelian groups.

\textbf{Remark (iii)}\qua If $H \le G \le \Gamma$ are three relatively hyperbolic
groups, such that  $G$ is fully quasi-convex in $\Gamma$, and $H$ is fully
quasi-convex in $G$, then $H$ is fully quasi-convex in $\Gamma$. Indeed, the
limit set of $H$ in $\Gamma$ is the image of the limit set of $H$ in $G$ by the
equivariant inclusion map $\partial (G) \hookrightarrow \partial(\G)$.

\begin{lemma}[`Full' intersection with parabolic subgroups]\label{lem;full}\Gap

  Let $\Gamma$ be a relatively hyperbolic
group with boundary $\partial \Gamma$, and $H$ be a fully quasi-convex subgroup.
Let $P$ be a parabolic subgroup of $\Gamma$. Then $P \cap H$ is either finite,
or of finite index in $P$.
\end{lemma}

Let $p \in \partial \Gamma$ the parabolic point fixed by $P$. Assume $P \cap H$ is not finite, so that
$p\in \Lambda H$. Then $p$ is in every translate of $\Lambda H$ by an element of $P$. The second point of
Definition \ref{def;full-q-c} shows that there are finitely many such translates: $P\cap H$ is of finite
index in $P$. \endproof

 \begin{prop}\label{prop;diam}
   Let $(\Gamma, \mathcal{G})$ be a relatively hyperbolic group, and $\partial \Gamma$ its
Bowditch boundary.
Let $H$ be
a quasi-convex subgroup of $\Gamma$, and $\Lambda H$ be its limit set in $\partial \Gamma$.
Let $(\gamma_n)_{n\in \mathbb{N}}$ be a sequence of elements of $\Gamma$ all in
distinct left cosets of $H$. Then there is a subsequence
$(\gamma_{\sigma(n)})$ such that $\gamma_{\sigma(n)}  \Lambda H$ uniformly
converges to a point.  \end{prop}

Unfortunately I do not know any purely dynamical proof of this proposition, that would only involve the
geometrically finite action on the boundary.

 There is a proper hyperbolic geodesic space $X$, with boundary $\partial \Gamma$,
on which $\Gamma$ acts properly discontinuously by isometries.
We assume that $\Lambda H$ contains two points $\xi_1$ and $\xi_2$, otherwise the result is a
consequence of the compactness of $\partial \Gamma$.
Let $B(\Lambda H)$ be the union of all the  bi-infinite geodesic between points of $\Lambda H$  in $X$,
and  $p$ be a point in it. Note that $B(\Lambda H)$ is quasi-convex in $X$, and
that $H$ acts on it properly discontinuously by isometries.  We prove that the
boundary  $\partial (B(\Lambda H))$ of $B(\Lambda H)$ is precisely $\Lambda H$.
Indeed, if $p_n$ is a sequence of points in $B(\Lambda H)$ going to infinity,
there are bi-infinite geodesics $(\xi_n,\zeta_n)$ containing each $p_i$, with
$\xi_n$ and $\zeta_n$ in $\Lambda H$. Let us extract a subsequence such that
$(\xi_n)_n$ converges to a point $\xi \in \partial (\G)$, and $\zeta_n \to \zeta
\in \partial (\G)$.  As $\Lambda H$ is closed, $\xi$  and $\zeta$ are in it, and
the sequence $(p_n)_n$ must converge to one of these two points (or both if they
are equal).

By our definition of quasi-convexity, $H$ acts on $\partial (B\Lambda H) = \Lambda H$ as a geometrically finite
convergence group.

To prove the proposition, it is enough to prove that a subsequence of
the sequence $\dist(\gamma_n^{-1} p, B( \Lambda H))$ tends to
infinity.  Indeed, by quasi-convexity of $B(\Lambda H)$ in $X$, for
all $\xi$ and $\zeta$ in $\Lambda H$, the Gromov products $(\gamma_n
\xi \cdot \gamma_n \zeta)_p$ are greater than $\dist(\gamma_n^{-1} p,
B( \Lambda H)) -K$, where $K$ depends only on $\delta$ and on the
quasi-convexity constant of $B(\Lambda H)$. Thus, we now want to prove
that a subsequence of $\dist(\gamma_n^{-1} p, B( \Lambda H))$ tends to
infinity.

 For all $n$, let $h_n \in H$ be such that $\dist(h_n p, \gamma_{n}^{-1} p)$ is
minimal among the distances $\dist(h p, \gamma_{n}^{-1} p)$, $h\in H$.
We prove the lemma:

\begin{lemma}
The sequence $(\dist(h_n p, \gamma_{n}^{-1} p))_n$ tends to infinity.
\end{lemma}

Indeed, if a
subsequence was bounded by a number $N$, then for infinitely many indexes, the
point $h_n^{-1} \gamma_{n}^{-1} p$ is in the ball of $X$ of center $p$ and of
radius $N$. Therefore, there exists $n$ and $m \neq n$ such that $h_n^{-1}
\gamma_{n}^{-1} = h_m^{-1} \gamma_{m}^{-1}$, which contradicts our hypothesis
that all the $\gamma_n$ are in distinct left cosets of $H$. \endproof

Let us resume the proof of Proposition \ref{prop;diam}.
For all $n$, let now $q_n$ be a point in $B(\Lambda H)$ such that
$\dist(\gamma_n^{-1} p, B( \Lambda H)) = \dist ( \gamma_n^{-1} p,  q_n)$.
 By the triangular inequality, $\dist (q_n, \gamma_n^{-1} p)
\geq \dist(h_n p, \gamma_{n}^{-1} p) - \dist (h_n p, q_n)$. If $(\dist
(h_n p, q_n))_n$ does not tend to infinity, then a subsequence of  $(\dist (q_n,
\gamma_n^{-1} p))_n$ tends to infinity and we are done.  Assume now that $(\dist
(h_n p, q_n))_n$ tends to infinity. After translation by $h_n^{-1}$, the
sequence $(\dist ( p, h_n^{-1} q_n))_n$ tends to infinity.
 Recall an usual result (Proposition 6.7 in
\cite{Brel}): given a $\Gamma$--invariant system of horofunctions
$(\rho_\xi)_{\xi\in \Pi}$, for the set $\Pi$ of bounded parabolic points in
$\partial \Gamma$, for all $t$,  there exists only finitely many horofunctions
$\rho_{\xi_1} \dots \rho_{\xi_k}$ such that $\rho_{\xi_i}(p) \geq t $.
As there are finitely many orbits of bounded parabolic points in $\Lambda H$, it
is possible to choose $t$ such that for every $\xi \in \Pi \cap \Lambda H$,
there exists $h\in H$ such that $\rho_{\xi}( h p) \geq t+1$.
 The group $H$,
 as a geometrically finite
 group,  acts co-compactly in the complement of a system of horoballs in
$B(\Lambda H)$ (Proposition 6.13 in \cite{Brel}). By definition of
the elements $h_n$, for all $h\in H$, one has  $\dist ( h p , h_n^{-1} q_n)
\geq \dist (p , h_n^{-1} q_n)$, and the latter tends to infinity.
Therefore the sequence $h_n^{-1} q_n$ leaves the complement of any system of
horoballs. In other words,  for all $M>0$, there exists $n_0$ such
that for all $n \ge n_0$, there is $i\in \{1, \dots, k\}$ such that
$\rho_{\xi_i} (h_n^{-1} q_n)\ge M$.

Therefore, one can extract a subsequence such that for some horofunction $\rho$
associated to a bounded parabolic point in $\Lambda H$, $\rho (h_n^{-1} q_n)$
tends to infinity. If $\dist (h_n^{-1} q_n, h_n^{-1} \gamma_n^{-1} p)$ remains
bounded, then $\rho (h_n^{-1} \gamma_n^{-1} p)$ tends to infinity, which is in
contradiction with Lemma 6.6 of \cite{Brel}, because $h_n^{-1} \gamma_n^{-1} p$
is in the $\Gamma$-orbit of $p$. Therefore a subsequence of $\dist (h_n^{-1}
q_n, h_n^{-1} \gamma_n^{-1} p)$ tends to infinity, and after translation by
$h_n$, this gives the result: a subsequence of $\dist (B(\Lambda H),
\gamma_n^{-1} p)$ tends to infinity. \endproof

The following statement appears in \cite{GrAsympt} and
also in \cite{Sh}, for hyperbolic groups. Note that this is no longer true for
(non fully) quasi-convex subgroups.

\begin{prop}[Intersection of fully quasi-convex
subgroups]\label{lem;interfullquasiconv}\Gap

 Let $\Gamma$ be a relatively hyperbolic group with boundary $\partial \Gamma$.
 If $H_1$ and $H_2$ are fully
quasi-convex subgroups of $\Gamma$, then $H_1 \cap H_2$ is fully quasi-convex,
moreover $\Lambda (H_1 \cap H_2) =\Lambda H_1 \cap \Lambda H_2$. \end{prop}

 As, for $i=1$ and $2$,  $H_i$ is a convergence group on $\Lambda H_i$, and as any sequence of distinct
 translates of $\Lambda H_i$ has empty intersection, the same is true for $H_1 \cap H_2$ on
 $\Lambda H_1 \cap \Lambda H_2$.

Let $p \in (\Lambda H_1 \cap \Lambda H_2)$ a parabolic point for $\Gamma$, and $P<\Gamma$ its stabilizer.
For $i=1$ and $2$, the group $H_i \cap P$ is maximal parabolic in $H_i$, hence infinite.
 By Lemma \ref{lem;full}, they are both of finite index in $P$, and therefore so is their intersection.
 Hence $p$ is a bounded parabolic point for $H_1 \cap H_2$ in $(\Lambda H_1 \cap \Lambda H_2)$.

Let $\xi \in (\Lambda H_1 \cap \Lambda H_2)$ be a conical limit point for
$\Gamma$. Then, by the first remark after the definition of quasi-convexity, it 
is a conical limit point for each of the $H_i$.

According to the definition of conical limit points, let $(\gamma_n)_{n\in
\mathbb{N}}$ be a sequence of elements in $\Gamma$ such that there exists
$\zeta$ and $\eta$ two distinct points in $\partial \Gamma$, with $\gamma_n \xi
\to \zeta$, and $\gamma_n \xi' \to \eta$ for all other $\xi'$.
There exists a subsequence of
$(\gamma_n)_{n\in \mathbb{N}}$ staying in a same left coset of $H_1$: if not, the
fact that two sequences $(\gamma_n \xi)_{n\in \mathbb{N}}$ and $(\gamma_n
\xi')_{n\in \mathbb{N}}$, for $\xi' \in \Lambda H_1 \setminus \{\xi\}$ converge to
two different points contradicts the Proposition \ref{prop;diam}.
By the same argument, there exists a subsequence of the previous
subsequence 
that remains in a
same left coset of $H_1$, and in a same left coset of $H_2$. Therefore it stays
in a same left coset of $H_1\cap H_2$; we can assume that we chose the sequence
$(\gamma_n)_{n\in \mathbb{N}}$ such that there exists $\gamma \in \Gamma$  and
$(h_n)_{n\in \mathbb{N}}$ a sequence of elements of $H_1 \cap H_2$, such that
$\forall n, \gamma_n = \gamma h_n$.

Therefore $h_n \xi \to \gamma^{-1}\zeta$, and $h_n \xi' \to \gamma^{-1}\eta$ for
all other $\xi'$. This means that $\xi \in \Lambda (H_1\cap H_2)$ is a conical
limit point for the action of $(H_1\cap H_2)$.   This ends the proof of
Proposition \ref{lem;interfullquasiconv}. \endproof

We emphasize the case of hyperbolic groups, studied by Bowditch in \cite{Bconf}.

\begin{prop}[Case of hyperbolic groups]\label{prop;hyp}\Gap

  In a hyperbolic group, a proper subgroup is quasi-convex in the sense of
quasi-convex subsets of a Cayley graph, if and only if it is fully
quasi-convex. \end{prop}

 B~Bowditch proved in \cite{Bconf} that a subgroup $H$ of a hyperbolic group 
$\Gamma$ is quasi-convex if and  only if it is hyperbolic with limit set
equivariantly homeomorphic to $\partial H$. It remains only to see that, if $H$
is quasi-convex in the classical sense, then the intersection of infinitely many
distinct translates $\bigcap_{n\in\mathbb{N}} (\gamma_n \partial H)$ is empty,
and we prove it by contradiction. Let us choose $\xi$ in
$\bigcap_{n\in\mathbb{N}} (\gamma_n \partial H)$. Then, there is $L>0$ depending
only on the quasi-convexity constant of $H$ in $\Gamma$, and there is, in each
coset $\gamma_n H$, an L-quasi-geodesic ray $r_n(t)$ tending to $\xi$. As they
converge to the same point in the boundary of a hyperbolic space, there is a
constant $D$ such that for all $i$ and $j$ we have: $ \exists t_{i,j} \,
\forall t>t_{i,j}, \, \exists t', \, {\rm dist}(r_i(t),r_j(t'))< D$. Let $N$ be
a number larger than the number of vertices in the a of radius $D$ in the Cayley
graph of $\Gamma$, and consider a point $r_1(T)$ with $T$ bigger than any
$t_{i,j}$, for $i,j \leq N$. Then each ray $r_i$, $i\leq N$,  has to pass
through the ball of radius $D$ centered in $r_1(T)$. By a pigeon hole argument,
we see that two of  them pass through the same vertex, but they were supposed to
be in disjoint cosets. \endproof

 Our point of view in Definition \ref{def;full-q-c} is a generalization of the
definitions in \cite{Bconf}, given for hyperbolic groups.

 \section{Boundary of an acylindrical graph of groups}

 Let $\Gamma$ be as the first or the second part of Theorem 0.1.
 We will say that we are in \emph{Case (1)} (resp.\ in \emph{Case (2)}) if
$\Gamma$ satisfies the first (resp.\ the second) assumption of Theorem 0.1.
However, we will need this distinction only for the proof of Proposition
\ref{coro;acyl}.

  Let $T$ be the Bass-Serre tree of the splitting, and $\tau$, a subtree of with
$T$ which is a fundamental   domain.
   We assume that the action of $\Gamma$ on $T$ is $k$--acylindrical for some $k \in \mathbb{N}$ (from Sela
   \cite{SelaAcyl}): the stabilizer of any segment of length $k$ is finite.

   We fix some notation: if $v$ is a vertex of $T$, $\Gamma_v$ is its
stabilizer in $\Gamma$. Similarly, for an edge $e$, we write $\Gamma_e$ for its
stabilizer. For a vertex $v$, $\G_v$ is relatively hyperbolic. This is by
assumption in Case (1), and in Case (2), if $\G_v$ is conjugated
to $A$, we consider that it is hyperbolic relative to itself; in this case the
space $\Sigma$ of Definition \ref{def;RHG} is just an horoball, and its Bowditch
boundary is a single point. For the existence of such a hyperbolic horoball, notice that 
the second definition of Bowditch \cite{Brel} indicates that the group $A*A$ is hyperbolic relative to the 
conjugates of both factors. Indeed we do not need to know the existence of such an horoball, but only that $A$ 
acts as a geometrically finite convergence group on a single point, which is trivial.

  \subsection{Definition of  $M$ as a set}

\textbf{Contribution of the vertices of $T$}

Let $\mathcal{V}(\tau)$ be the set of vertices of $\tau$.
For a vertex $v\in \mathcal{V}(\tau)$, the group $\Gamma_v$
is by assumption a relatively hyperbolic group and we denote by
$\partial(\Gamma_v)$ a compact space homeomorphic to its Bowditch boundary. Thus,
$\G_v$ is a geometrically finite convergence
group on $\partial(\G_v)$.

 We
set $ \Omega $ to be
$\Gamma \times \left(\bigsqcup_{v\in \mathcal{V}(\tau)}
\partial(\G_v)\right) $  divided by the natural relation
$$(\gamma_1,
x_1)=(\gamma_2,x_2) \; \hbox{ if } \; \exists v \in \mathcal{V}(\tau), \;
x_i \in \partial (\G_v), \gamma_2^{-1}\gamma_1 \in \G_v,
\gamma_2^{-1}\gamma_1 x_1=x_2.$$

In particular, for each $v \in \tau$, the compact space $\partial \G_v$ naturally
embeds in $\Omega$ as the image of $\{1\} \times \partial \G_v$. We identify it
with its image. The group $\Gamma$ naturally acts on the left on $\Omega$.
The compact space $\gamma \partial(\G_v)$ is called the boundary of the vertex
stabilizer $\Gamma_{\gamma v}$.

\rk{Contribution of the edges of $T$}

Each edge will allow us to glue together boundaries of vertex stabilizers along 
the limit sets of the stabilizer of the edge. We explain precisely how.

 For an edge $e=(v_1,v_2)$ in $\tau$, the group $\G_e$  embeds
as a quasi-convex subgroup in both  $\G_{v_i}$, $i=1,2$.  Thus, by
definition of quasi-convexity, these embedings define equivariant maps
$\Lambda_{(e,v_i)}: \partial (\G_e) \hookrightarrow \partial(\G_{v_i})$, where
$\partial(\G_e)$ is the Bowditch boundary of the relatively hyperbolic group
$\G_e$. Similar maps are defined by translation, for edges in $T \setminus
\tau$.

The equivalence
relation $\sim$ on $\Omega$ is the transitive closure of the following: for
$v$ and $v'$ are vertices of $T$, the points $\xi \in \partial(\G_v)$  and $\xi'
\in \partial(\G_{v'})$ are equivalent in  $\Omega$ if there is an edge $e$
between $v$ and $v'$, and a point $x \in \partial(\G_e)$ satisfying
simultaneously $\xi = \Lambda_{(e,v)}(x)$ and $\xi'= \Lambda_{(e,v')}(x)$.

\begin{lemma}

Let $\pi$ be the projection corresponding to the quotient:  $\pi: \Omega \to \Omega/_{\sim} $.
 For all vertex $v$, the restriction of $\pi$ on $\partial(\G_v)$ is injective.

\end{lemma}

Let $\xi$ and $\xi'$  be two points of $\Omega$, both of them
being in the boundary of a vertex stabilizer $\partial (\G_v)$. If they are
 equivalent for the relation above, then there is a sequence of consecutive
 edges $e_1=(v,v_1), e_2= (v_1,v_2) \dots e_n = (v_{n-1}, v)$, the first one
 starting at $v_0=v$ and the last one ending at $v_n=v$, and a sequence of points
 $\xi_i \in \partial (\G_{v_i})$, for $i\leq n-1$,  such that, for all $i$,
 there exists $x_i \in \partial(\G_{e_i})$, satisfying
 $\xi_i = \Lambda_{(e_i,v_{i-1})}(x_i)$ and
 $\xi_{i+1}= \Lambda_{(e_i,v_i)}(x_i)$. As $T$ is a tree, it contains no simple
 loop, and there exists an index $i$ such that $v_{i-1} = v_{i+1}$. As, for all
 $j$, the maps $\Lambda_{(e_j,v_j)}$ are injective, the points $\xi_{i-1}$ and
 $\xi_{i+1}$ are the same in $\partial (\G_(v_{i-1}))$, and inductively, we see
 that $\xi$ and $\xi'$ are the same point. This proves the claim. \endproof

 Note that the group $\Gamma$
acts on the left on $\Omega/_{\sim}$.
Let $\partial T$ be the (visual) boundary of the tree $T$: it is the space of
the rays in $T$ starting at a given base point; let us recall that for its
topology, a sequence of rays $(\rho_n)$ converges to a given ray $\rho$, if
$\rho_n$ and $\rho$ share arbitrarily large prefixes, for $n$ large enough. We
define $M$ as a set: $$   M= \partial T \sqcup (\Omega/_{\sim}). $$

As before, let $\pi$ be the projection corresponding to the quotient:  $\pi: \Omega \to \Omega/_{\sim} $.
For a given edge $e$ with vertices $v_1$ and $v_2$, the two maps
$\pi \circ \Lambda_{(e,v_i)}: \partial(\G_e) \to \Omega/_{\sim}$ $(i=1,2)$,
are two equal homeomorphisms on their common image. We identify this image with
the Bowditch boundary of $\G_e$, $\partial(\G_e)$, and we call
this compact space, the boundary of the edge stabilizer  $\G_e$.

  \subsection{Domains}

Let $\mathcal{V}(T)$ be the set of vertices of $T$.
 We still denote by $\pi$ the projection:  $\pi: \Omega \to \Omega/_{\sim} $.
 Let $\xi \in \Omega/_{\sim}$. We define the \emph{domain} of $\xi$, to be
 $ D(\xi) \; = \; \{v\in \mathcal{V}(T) \: | \: \xi \in \pi(\partial(\G_v))
\}$. As we want uniform notations  for all points in $M$, we say that the
\emph{domain of a point} $\xi \in \partial T$ is $\{\xi\}$ itself.

 \begin{prop}[Domains are bounded]\label{coro;acyl}\Gap

        For all $\xi \in \Omega/_{\sim}$, $D(\xi)$ is convex in $T$, its diameter is bounded by the
        acylindricity constant, and  the intersection of two distinct domains is finite. The quotient of $D(\xi)$
        by the stabilizer of $\xi$ is finite.
 \end{prop}

\begin{rem} In Case (1), we will even prove that domains are \emph{finite}, but
this is false in Case (2).
\end{rem}

  The equivalence $\sim$ in $\Omega$ is the transitive closure of a relation involving points in boundaries
  of adjacent vertices, hence domains are convex.

 \rk{End of the proof in Case (2)}  As $P$ is a maximal parabolic subgroup of
$G$,  its limit set is a single point: $\partial(P)$ is one point belonging to
the boundary of only one stabilizer of an edge adjacent to the vertex
$v_G$ stabilized by $G$. Hence, the domain of $\xi = \partial (\G_{v_A})$ is
$\{v_A\}\cup Link(v_A)$, that is $v_A$ with all its neighbours, whereas the
domain of a point $\zeta$  which is not a translate of $\partial (\G_{v_A})$, is
only one single vertex.

 Domains have therefore diameter bounded by $2$, and any two of them intersect
only on one point. For the last assertion, note that $A$  stabilizes the point
$\partial(\G_{v_A})$, and acts transitively on the edges adjacent to $v_A$. This
proves the lemma in {Case (2)}.

In {Case (1)}, we need a lemma:

 \begin{lemma} \label{lem;domains}
In {Case (1)}, let $\xi \in \Omega/_{\sim}$. The stabilizer of any finite subtree of $D(\xi)$ is
infinite.

If a subtree, whose vertices are $\{v_1,\dots, v_n\}$, is in $D(\xi)$, then
there exists a group $H$ embedded in each of the vertex
stabilizers $\G_{v_i}$ as a fully quasi-convex subgroup, with $\xi$ in
its limit set.
 \end{lemma}

 The first assertion is clearly a
consequence of the second one, we will prove the latter by induction.

 If $n=1$, $H$ is the vertex stabilizer.
 For larger $n$, re-index the vertices so that $v_n$ is a final leaf of the subtree $\{v_1,\dots,v_n\}$,
with unique neighbor $v_{n-1}$. Let $e$ be the edge $\{v_{n-1},v_n\}$.
 The induction gives $H_{n-1}$, a subgroup of the stabilizers of each $v_i$, $i\leq n-1$, and with
 $\xi \in \partial H_{n-1}$. As $\xi \in \partial(\G_{v_n})$, it is in
$\partial(\G_e)$, and we have  $\xi \in \partial H_{n-1} \cap \partial(\G_e)$.
 By Proposition \ref{lem;interfullquasiconv}, $H_{n-1}\cap \G_e$ is a fully
quasi-convex subgroup of $\G_{v_{n-1}}$, and therefore, it is a a fully
quasi-convex subgroup of $\G_{e}$, and of $H_{n-1}$. Therefore, (see Remark
(iii)), it is a fully quasi-convex subgroup of $\G_{v_n}$, and of each
of the $\G_{i}$, for $i\leq (n-1)$,  with $\xi$ in its limit set. It is then
adequate for $H$; this proves the claim, and Lemma \ref{lem;domains}. \endproof

\proof[{End of the proof of Prop. \ref{coro;acyl} in Case (1)}]
By Lemma \ref{lem;domains}, each segment in $D(\xi)$ has an infinite stabilizer,
hence by $k$--acylindricity, $Diam(D(\xi))\leq k $.   Domains
 are bounded, and they are locally finite because of the second
requirement of Definition  \ref{def;full-q-c}, therefore they are finite.
The other assertions are now obvious. \endproof

   \subsection{Definition of neighborhoods in $M$}

 We will describe $(W_n(\xi))_{n\in \mathbb{N}, \xi \in M}$, a family of subsets
of $M$, and prove that it  generates an topology (Theorem
\ref{thm;topo}) which is suitable for our purpose.

For a vertex $v$, and an open subset $U$ of $\partial(\G_{v})$, let $T_{v,U}$ be
the subtree whose vertices $w$ are such that $[v,w]$ starts by an edge $e$ with
$\partial(\G_e) \cap U \neq \emptyset$.

For each vertex $v$ in $T$, let us choose $\mathcal{U}(v)$, a countable basis of
open neighborhoods of $\partial(\G_v)$, seen as the Bowditch boundary of $\G_v$.
Without loss of generality, we can assume that for all $v$, the
collection of open subsets $\mathcal{U}(v)$ contains $\partial(\G_v)$ itself.

Let $\xi$ be in $\Omega/_{\sim}$, and  $D(\xi)=\{v_1, \dots, v_n, \dots \} =
(v_i)_{i\in I}$. Here, the set $I$ is a subset of $\mathbb{N}$. For each $i\in
I$, let $U_i \subset \partial(\G_{v_i})$ be an element of $\mathcal{U}(v_i)$,
containing $\xi$, such that for all but finitely many indices  $i\in I$,  
$U_i=\partial(\G_{v_i})$.

 The set $W_{(U_i)_{i\in I}}(\xi)$ is the disjoint union of three subsets
$ W_{(U_i)_{i\in I}}(\xi) = A \cup B \cup C$:

\medskip

 $\bullet$\qua $A =  \bigcap _{i\in I}  \partial (T_{v_i,U_i})$,

 $\bullet$\qua  $ B = \{ \zeta \in (\Omega/_{\sim}) \setminus  (\bigcup_{i\in I}
\partial(\G_{v_i})) \: | \;  D(\zeta) \subset
\bigcap_{i\in I} T_{v_i,U_i} \}$

 $\bullet$\qua  $C=  \{\zeta \in \bigcup_{j\in I} \partial(\G_{v_j}) \, | \;  \zeta \in \bigcap_{m\in I |\zeta \in \partial(\G_{v_m})}  U_m \}$.

\medskip

{\bf Remark}\qua The set of elements of $\Omega/_{\sim}$ is not countable,
nevertheless, the set of different possible domains is countable. Indeed a
domain is a finite subset of vertices of $T$ or the star of a vertex of $T$, and
this makes only countably many possibilities. The set $W_{(U_i)_{i\in I}}(\xi)$
is completely defined by the data of the domain of $\xi$, the data of a finite
subset $J$ of $I$, and the data of an element of $\mathcal{U}(v_j)$ for each
index $j\in J$. Therefore, there are only countably many different sets
$W_{(U_i)_{i\in I}}(\xi)$, for $\xi \in \Omega/_{\sim}$, and $U_i \in
\mathcal{U}(v_i)$, $v_i \in D(\xi)$.  For each $\xi$ we choose an arbitrary
order  and denote them $W_m(\xi)$.

 Let us choose $v_0$ a base point in $T$. For $\xi \in \partial T$, we define
 the subtree $T_m(\xi)$: it consists of the vertices $w$ such that
$[v_0,w] \cap [v_0,\xi)$ has length  bigger than $m$. We set
 $W_m(\xi) =  \{\zeta \in M \; |\; D(\zeta)\subset T_m(\xi) \}$. Up to a
shift in the indexes, this does not depend on $v_0$, for $m$ large enough.

\begin{lemma}[Avoiding an edge]\label{lem;filtrage}\Gap

  Let $\xi$ be a point in $M$, and $e$ an edge in $T$ with at least one vertex not in $D(\xi)$. Then,
  there exists an integer $n$ such that $W_n(\xi) \cap \partial(\G_e) =
\emptyset$. \end{lemma}

If $\xi$ is in $\partial T$ the claim is obvious.
 If $\xi \in \Omega/_{\sim}$, as $T$ is  a tree, there is a unique segment from
the convex $D(\xi)$ to $e$.  Let $v$ be the vertex of $D(\xi)$ where this path
starts, and $e_0$ be its first edge. It is enough  to find a neighborhood  of
$\xi$ in $\partial(\G_v)$ that misses $\partial(\G_{e_0})$. As one vertex of
$e_0$ is not in $D(\xi)$, $\xi$ is not in $\partial(\G_{e_0})$, which is
compact. Hence such a separating neighborhood exists. \endproof

   \subsection{Topology of $M$}

  In the following, we consider the smallest topology $\mathcal{T}$ on $M$ such that the family of sets
$\{W_n(\xi); \, n\in \mathbb{N}, \, \xi \in M \}$, with the notations above, are open subsets of $M$.

  \begin{lemma}\label{lem;haus}
    The topology $\mathcal{T}$ is Hausdorff.
  \end{lemma}

 Let $\xi$ and $\zeta$ two points in $M$.
If the subtrees $D(\xi)$ and $D(\zeta)$ are disjoint, there is an edge $e$ that
separates them in $T$, and Lemma \ref{lem;filtrage} gives two neighborhoods of
the points that do not intersect.  Even if $D(\xi)\cap  D(\zeta)$ is non-empty, it is
nonetheless finite (Proposition \ref{coro;acyl}). In each of its vertex $v_i$, we can choose
disjoint neighborhoods $U_i$ and $V_i$ for the two points. This gives rise to sets
$W_n(\xi)$ and $W_m(\zeta)$ which are separated. \endproof

  \begin{lemma}[Filtration]\label{lem;filtration}\Gap

   For every $\xi \in M$, every integer $n$, and every $\zeta \in W_n(\xi)$,
there exists $m$ such that $W_m(\zeta) \subset W_n(\xi)$.
 \end{lemma}

If $D(\zeta)$ and $D(\xi)$ are disjoint, again, Lemma \ref{lem;filtrage} gives
a neighborhood of $\zeta$, $W_m(\zeta)$ that do not meet $\partial(\G_e)$, whereas
$\partial(\G_e) \subset W_n(\xi)$, because $\zeta \in W_n(\xi)$. By definition of
our family of neighborhoods, $W_m(\zeta)\subset W_n(\xi)$.

If the domains of $\xi$ and $\zeta$ have a non-trivial intersection, either the
two points are equal (and there is nothing to prove), or the intersection is
finite (Prop. \ref{coro;acyl}). Let $(v_i)_{i\in I}=D(\xi)$,  let
$(U_i)_{i\in I}$ be such that $W_n(\xi) = W_{(U_i)_i}(\xi)$, and let $J\subset
I$ be such that $D(\xi)\cap D(\zeta) = (v_j)_{j\in J}$.  In this case, we can
choose, for all $j\in J$, a neighbourhood of $\zeta$ in $\partial(\G_j)$,
$U'_j \subset U_j$ such that $U'_j$ do not meet the boundary of the stabilizer
of an edge $(v_j,v_i)$ for any $i\in I\subset J$; this gives $W_m(\zeta)\subset
W_n(\xi)$. \endproof

\begin{coro}\label{coro;system}
The family  $\{ W_n(\xi)\}_{ n\in \mathbb{N}, \xi \in M }$ is a fundamental
system of open neighborhoods of $M$ for the topology  $\mathcal{T}$.
\end{coro}

 It is enough to show that the intersection of two such sets is equal to the
union of some other ones. Let $W_{n_1}(\xi_1)$ and $W_{n_2}(\xi_2)$ be in the
family. Let $\zeta$ be in their intersection. Lemma \ref{lem;filtration} gives
$W_{(U_j)_j}(\zeta)\subset W_{n_1}(\xi_1)$ and  $W_{(V_j)_j}(\zeta)\subset
W_{n_2}(\xi_2)$. As $W_{(U_j)_j}(\zeta) \cap W_{(V_j)_j}(\zeta) = W_{(U_j\cap
V_j)_j}(\zeta)$, we get an integer $m_{\zeta}$ such that $W_{m_{\zeta}}(\zeta)$
is included in both $W_{n_i}(\xi_i)$. Therefore, $W_{n_1}(\xi_1) \cap
W_{n_2}(\xi_2) = \bigcup_{\zeta \in W_{n_1}(\xi_1) \cap W_{n_2} (\xi_2)}
W_{m_{\zeta}} (\zeta)$. \endproof

\begin{coro}
Recall that $\pi$ be the projection corresponding to the quotient:  $\pi:
\Omega \to \Omega/_{\sim} $.  For all vertex $v$, the restriction of $\pi$ on
$\partial(\G_v)$ is continuous.
\end{coro}

Let $\xi$ be in $\partial(\G_v)$, and let $(\xi_n)_n$ be a sequence of elements
of $\partial(\G_v)$ converging to $\xi$ for the topology of $\partial(\G_v)$.
Let $(U^n)_n$ be a system of neighbourhoods of $\xi$ in $\partial(\G_v)$, such
that for all $n$, for all $n'\geq n$, $\xi_{n'} \in U^n$.
 Let $D(\pi(\xi))= \{v, v_2,
\dots \} $ in $T$, and consider $W_m=W_{(U_i(m))} (\pi(\xi))$, such that
$U_1(m)\subset U^n$. By definition, $W_{(U_i(m))} (\pi(\xi))\cap
\pi(\partial(\G_v))$ is the image by  $\pi$ of an open subset of $U_1(m)$
containing $\xi$. Therefore, by property of fundamrental systems of
neighbourhoods, $\pi(\xi_n)$ converges to $\pi(\xi)$. Therefore $\pi$ is
continuous. \endproof

From now, we identify $\xi$ and $\pi(\xi)$ in such situation.

  \begin{lemma}\label{lem;reg}
    The topology  $\mathcal{T}$ is regular, that is, for all $\xi$, for all $m$, there exists
$n$ such that $\overline{U_n(\xi)}\subset U_m(\xi)$.
 \end{lemma}

 In the case of $\xi \in  \partial T$, the closure of $W_n(\xi)$  is contained
in $W'_n(\xi) =  \{\zeta \in M | D(\zeta)\cap T_n(\xi) \neq \emptyset \}$
(compare with the definition of $W_n(\xi)$). As, by Proposition \ref{coro;acyl},
domains have uniformly bounded diameters, we see that for arbitrary $m$, if $n$
is large enough, $\overline{W_n(\xi)} \subset W_m(\xi)$.

 In the case of $\xi \in \Omega/_{\sim}$, $\overline{W_{(U_i)_i}(\xi)}\setminus
W_{(U_i)_i}(\xi) $ contains only points in the boundaries of vertices of
$D(\xi)$, and those are in the closure of the $U_i$ (which is non-empty only for
finitely many $i$), and in the boundary (not in $U_i$) of edges meeting $U_i
\setminus \{\xi\}$. Therefore, given $V_i\subset \partial(\G_{v_i})$,  with strict
inclusion only for finitely many indices, if we choose the $U_i$ small enough to
miss the boundary of every edge non contained in $V_i$, except the ones meeting
$\xi$ itself, we have $\overline{W_{(U_i)_i}(\xi)} \subset W_{(V_i)_i}(\xi)$.
\endproof

 \begin{theo}\label{thm;topo}
    Let $\Gamma$ be as in Theorem 0.1.
   With the notations above, $\{W_n(\xi); \, n\in \mathbb{N}, \, \xi \in M \}$
 is a base of a topology that makes $M$ a perfect
metrizable compact space, with the following convergence criterion:
 $ \left(    \xi_n \to \xi  \right)  \Longleftrightarrow \left(
\forall n \exists m_0 \forall m>m_0, \xi_m \in W_n(\xi) \right)$.  \end{theo}

 The topology is, by construction, second countable, separable. As it is
also Hausdorff (Lemma \ref{lem;haus})  and regular (Lemma \ref{lem;reg}), it is metrizable.
The convergence criterion is an
immediate consequence of Corollary \ref{coro;system}. Let us prove that it is
sequentially compact. Let $(\xi_n)_{n\in \mathbb{N}}$ be a sequence in $M$, we
want to extract a converging subsequence. Let us choose $v$ a vertex in $T$, and
for every $n$, $v_n \in D(\xi_n)$ minimizing the distance to $v$ (if $\xi_n \in
\partial T$, then $v_n=\xi_n$). There are two possibilities (up to extracting a
subsequence): either the Gromov products $(v_n\cdot v_m)_v$ remain bounded, or
they go to infinity. In the second case, the sequence $(v_n)_n$ converges to  a
point in $\partial T$, and by our convergence criterion, we see that $(\xi_n)_n$
converges to this point (seen in $\partial T \subset M$). In the first case,
after extraction of a subsequence, one can assume that the Gromov products
$(v_n\cdot v_m)_v$ are  constant equal to a number $N$. Let $g_n$ be a geodesic
segment or  a geodesic ray between $v$ and $v_n$. there is a segment $g=[v,v']$
of length $N$, which is contained in every $g_n$, and for all distinct $n$ and
$m$, $g_n$ and $g_m$ do not have a prefix longer than $g$.

 Either there is a subsequence so that $g_{n_k} = g$ for
all $n_k$, and as $\partial \G_{v'}$ is compact, this gives the result, or there
is a subsequence such that every $g_{n_k}$ is strictly longer than $g$. Let
$e_{n_k}$ be the edge of $g_{n_k}$ following $v'$. All the
$e_{n_k}$ are distinct, therefore, by Proposition \ref{prop;diam}, one can extract another subsequence
such that  the
sequence of the boundaries of their stabilizers converge to a single point of $\partial
\G_{v'}$. The convergence criterion indicates that the subsequence of $(\xi_{n_k})_n$
converges to this point.

Therefore, $M$ is sequentially compact and metrisable, hence it is compact. It
is perfect since $\partial T$ has no isolated point, and accumulates everywhere.
\endproof

\begin{theo}[Topological dimension of $M$]
\label{thm;dim}{\rm [Theorem 0.2]}

  $dim(M) \leq \max_{v,e} \{dim(\partial(\G_v)), dim(\partial(\G_e))+1 \}$.
\end{theo}

It is enough to show that every point has arbitrarily small neighborhoods whose
boundaries have topological dimension at most $(n-1)$ (see the book
\cite{HW}, where this property is set as a definition).

If $\xi \in \partial T$, the closure of $W_n(\xi)$  is contained
in $W'_n(\xi) =  \{\zeta \in M | D(\zeta)\cap T_n(\xi) \neq \emptyset \}$
(compare with the definition of $W_n(\xi)$). The boundary
of $W_n(\xi)$ is therefore a compact subspace of the boundary of the stabilizer of
the unique  edge that has one and only one vertex in $T_n(\xi)$; the boundary
of $W_n(\xi)$ has dimension at most $\max_{e} \{dim(\partial(\G_e))\}$.

If $\xi \in \Omega /_{\sim} $,  $\overline{W_{(U_i)_i}(\xi)}\setminus
W_{(U_i)_i}(\xi) $ contains only points in the boundaries of vertices of
$D(\xi)$, and those are in the closure of the $U_i$ (which is non-empty only for
finitely many $i$), and in the boundaries (not in $U_i$) of stabilizers of edges
that meet $U_i \setminus \{\xi\}$. Hence, the boundary of a neighborhood
$W_n(\xi)$ is the union of boundaries of neighborhoods of $\xi$ in
$\partial(\G_{v_i})$ and of a compact subspace of the boundary of countably many
stabilizers of edges. As the dimension of a countable union of compact spaces of
dimension at most $n$ is of dimension at most $n$ (Theorem III.2 in \cite{HW}),
its dimension is therefore at most $\max_{v,e} \{dim(\partial(\G_v))-1,
dim(\partial(\G_e))\}$. This proves the claim. \endproof

\section{Dynamic of $\Gamma$ on $M$}

We assume the same hypothesis as for Theorem  \ref{thm;topo}. We first prove
two lemmas, and then we prove the different assertions of Theorem
\ref{theo;dynamic}.

\begin{lemma}[Large translations]\label{lem;gdetrans}\Gap

  Let $(\gamma_n)_{n\in \mathbb{N}}$ be a sequence in $\Gamma$. Assume that, for
some (hence any) vertex $v_0\in T$,  ${\rm dist}(v_0, \gamma_n v_0) \to \infty$. 
Then, there is a subsequence $(\gamma_{\sigma(n)})_{n\in \mathbb{N}}$, there is
a point $\zeta \in M$, and a point $\zeta' \in \partial T$, such that for all
compact subspace $K \subset (M\setminus\{\zeta'\})$,
one has $\gamma_{\sigma(n)}K \to \zeta$ uniformly.
\end{lemma}

Let $\xi_0$ be in $\partial(\G_{v_0})$.
Using the sequential compactness of $M$, we choose a
subsequence $(\gamma_{\sigma(n)})_{n\in \mathbb{N}}$ such that
$(\gamma_{\sigma(n)}\xi_0)_n$ converges to a point $\zeta$ in $M$; we still
have ${\rm dist}(v_0, \gamma_{\sigma(n)} v_0) \to \infty$.

Let $v_1$ be another
vertex in $T$. The lengths of the segments $[\gamma_n v_0, \gamma_n v_1]$ are all equal
to the length of $[v_0, v_1]$, therefore, for all $m$, there is $n_m$ such that for all $n>n_m$,
the segments $[v_0, \gamma_{\sigma(n)} v_0]$ and
$[v_0,\gamma_{\sigma(n)} v_1]$  have a common prefix of length more than $m$.

Let $\zeta_1, \zeta_2$ in $\partial T$. The center of the triangle $(v_0, \zeta_1,
\zeta_2)$ is a vertex $v$ in $T$. Therefore, for all $m\ge 0$, the segments $[v_0,
\gamma_{\sigma(n)} v_0]$ and $[v_0,\gamma_{\sigma(n)} v]$ coincide on a
subsegment of length more than $m$, for sufficiently large integers $n$. This means that for
at least one of the $\zeta_i$, the ray  $[v_0, \gamma_{\sigma(n)} \zeta_i)$
has a common prefix with  $[v_0,\gamma_{\sigma(n)} v_0]$ of length at least
$m$. By convergence criterion, $(\gamma_{\sigma(n)} \zeta_i)$ converges to
$\zeta$. Therefore there exits $\zeta'$ in $\partial \Gamma$, such that any other point
$\zeta" \in (\partial T \setminus \{\zeta'\})$, satisfies  $\gamma_{\sigma(n)}\zeta" \to \zeta$.

Let $K$ be a compact subspace of  $(M\setminus\{\zeta'\})$. There exists a vertex
$v_0$, a point $\xi \in \partial T$, and a neighborhood $W_{m}(\xi)$ (see the
definition in the section above, where $v_0$ is the base point) of $\xi$
containing $K$, not containing $\zeta'$. Let $v$ be on the ray $[v_0,\xi)$, at
distance $m$ from $v_0$. Then for all points $\xi'$ in $W_{m}(\xi)$ the ray
$[v_0, \xi')$ has the prefix $[v_0, v]$.
 As the sequence $(\gamma_{\sigma(n)} \partial \G_v)_{n\in \mathbb{N}}$
 uniformly converges to $\zeta$, the sequence   $(\gamma_{\sigma(n)}  W_m(\xi))_{n\in
\mathbb{N}} $ uniformly converges to this point. Therefore, the convergence is
uniform on $K$. \endproof

\begin{lemma}[Small translation]\label{lem;petitetrans}\Gap

 Let $(\gamma_n)_{n\in \mathbb{N}}$ be a sequence of distinct elements of
$\Gamma$, and assume that for some (hence any) vertex $v_0$, the sequence
$(\gamma_n v_0)_n$ is bounded in $T$. Then there exists a subsequence
$(\gamma_{\sigma(n)})_{n\in \mathbb{N}}$, a vertex $v$, a point $\zeta \in
\partial(\G_v)$, and another point $\zeta' \in \Omega/_{\sim}$, such that, for
all compact subspace $K$ of $M\setminus\{\zeta'\}$, one has $\gamma_{\sigma(n)} K
\to \zeta$ uniformly. \end{lemma}

We distinguish two cases.
First, we assume that for some vertex $v$, and for some element $\gamma \in
 \Gamma$, there exists a subsequence such that $\gamma_n
= h_n \gamma$, with $h_n \in \G_v$ for all $n$. In such a case, we can extract
again a subsequence (but, without loss of generality, we still denote it by
$(\gamma_n)_n$) such that there exists a point $\zeta' \in \partial
(\G_{\gamma^{-1}v})$ and a point $\zeta \in \partial(\G_v)$, such that for every
compact subspace $K_{\gamma^{-1}v} \subset \partial (\G_{\gamma^{-1}v}) \setminus
\{\zeta'\}$, our subsequence of $\gamma_n K_{\gamma^{-1}v}$ converges to $\zeta$
uniformly.

Assume that $\zeta'$ is not a parabolic point for $\G_v$ in $\partial(\G_v)$.
 For any vertex $w$ in $D(\gamma \zeta')$, let $e$ be the first edge
of the segment $[v,w]$. The boundary of its stabilizer contains $\zeta'$. The
elements $h_n$ are all, except finitely many, in the same left coset of
$Stab(e)$, otherwise, as $h_n \gamma \zeta'$ and $h_n \xi$ go to different
points, for all $\xi\neq \gamma \zeta'$ in $\partial(\G_e) \setminus \{\zeta'\}$
(which is non empty since $\zeta'$ is not parabolic), we get a contradiction
with Proposition \ref{prop;diam}. Therefore, we can extract a subsequence (but,
without loss of generality, we still denote it by $(\gamma_n)_n$) such that,
for each vertex $\gamma^{-1}w \in D(\zeta')$, for each compact
subspace $K_{\gamma^{-1} w}$ of $\partial (\G_{\gamma^{-1}w})$, not containing
$\zeta'$, the sequence $\gamma_{n} K_{\gamma^{-1}w}$ converges  to $\zeta$
uniformly.
Assume now that  $\zeta'$ is a parabolic point for $\G_v$ in $\partial(\G_v)$.
Then $h_n(\gamma \zeta')$ do converge to $\zeta$, otherwise, $\zeta'$ would
be a conical limit point. Therefore, for all vertex $\gamma^{-1}w \in
D(\zeta')\setminus \{\gamma^{-1}v\}$, the sequence $\gamma_{n}
\partial(\G_{\gamma^{-1}w})$ converges to $\zeta$ uniformly.

Therefore, if $v'$ is a vertex
not in the domain of $\zeta'$,  the path from $ \gamma^{-1}v$ to $v'$ contains
an edge such that the boundary  of its stabilizer is a compact space
$K_{\gamma^{-1}w}$ satisfying: $\gamma_{n} K_{\gamma^{-1}w} \to \zeta'$
uniformly. Let $K$ be a compact subspace of $M\setminus\{\zeta'\}$. For each $v_i
\in D(\zeta')$, there exists a compact space $K_i\subset \partial(\G_{v_i})
\setminus\{\zeta'\}$, $K\cap \partial(\G_{v_i}) \subset K_i$ such that for all
other point $\xi$ of $K$, the unique ray in $T$ from $D(\zeta')$ that converges
to $\xi$ contains an edge such that the boundary of its stabilizer is contained
in some $K_i$. Therefore, $\gamma_{n} K \to \zeta'$ uniformly.

We turn now to the second case, where such a subsequence does not exists.
Nevertheless, after extraction, we can assume that the distance
$\dist (v_0, \gamma_n v_0)$ is constant.
Let $v$ be the vertex such that there exists a subsequence
$(\gamma_{\sigma(n)})_{n\in \mathbb{N}}$ with the property that some segments
 $[v_0, \gamma_{\sigma(n)} v_0]$ have a common prefix $[v_0, v]$, and
the edges $e_{\sigma (n)} \subset [v_0, \gamma_{\sigma(n)} v_0]$ located just
after $v$, are all distinct. By Proposition \ref{prop;diam}, one can extract
 a subsequence $(e_{\sigma'(n)})_n $ such that the boundaries of the
stabilizers of these edges converge to some point $\zeta \in \partial(\G_v)$.
By our convergence criterion,
$\gamma_{\sigma'(n)} \partial(\G_{v_0})$ uniformly converges  to $\zeta$.

Let $\xi$ be a point in $\partial T$.  We claim that $v$ is not in the
ray  $[\gamma_{\sigma'(n)} v_0,\gamma_{\sigma'(n)} \xi)$, for $n$ sufficiently
large. If it was, there would be a subsequence satisfying:
$\gamma_{\sigma'(n)}^{-1} v$ is constant on a vertex $w$ of the ray $ [v_0,
\xi)$, that is, $\gamma_{\sigma'(n)}^{-1} =  h_n \gamma $, where $h_n \in \G_{w}
$. Therefore, $\gamma_{\sigma'(n)}w$ equals $v$ for all $n$. In other
words, for all $n$ there exists $h_n$ in $\G_w$ such that
$\gamma_{\sigma'(n)} = h_n \gamma_{\sigma'(0)}$.  This contradicts our
assumption that we are not in the first case, and this proves the claim.

If $d = \dist (\gamma_{\sigma' (n)} v_0, v)$ (which is constant by assumption),
we choose the neighborhood of $\xi$ defined by $W_{d+1}(\xi)$ (here $v_0$ is the
base point). Then, for each point in  $\gamma_{\sigma' (n)} W_{d+1}(\xi)$, the
unique path in $T$ from $v_0$ to this point contains $e_n$.
Therefore, $\gamma_{\sigma' (n)} W_{d+1}(\xi)$, uniformly converges to $\zeta$.

Let $\xi$ be now a point in the boundary of the stabilizer of a vertex $v'$.
Again, for the same reason, the vertex $v$ is not in $[\gamma_{\sigma'(n)}
v_0,\gamma'_{\sigma(n)}   v']$ for $n$ large enough.
Therefore the unique path from $v$ to $\gamma_{\sigma'(n)}   v'$
contains the edge $e_{\sigma'(n)}$. If $\gamma_{\sigma'(n)} \xi$ is not in $
\partial (\G_{e_{\sigma'(n)}})$, for all $n$ sufficiently large, then there
exists a neighborhood $N$ of $\xi$ such that the convergence $\gamma_{\sigma'
(n)} N \to \zeta$ is uniform.   If $\gamma_{\sigma'(n)} \xi$ is in $ \partial
(\G_{e_{\sigma'(n)}})$, then there exists another vertex $v_n''$ in $D(\xi)$
such that $\gamma_{\sigma'(n)}(v_n'') = v$. If $D(\xi)$ is finite, after
extracting another subsequence, we see that we are in the first case, but we
supposed we were not. If $D(\xi)$ is infinite, we are in case (2) of the main
theorem, and $D(\xi)$ is exactly the star of a vertex $v"$. If $v$ is in the
orbit of the vertex stabilized by the group $A$, again, necessarily
$\gamma_{\sigma'(n)}(v") = v$. If $v$ is not in this orbit,
$\gamma_{\sigma'(n)}^{-1} v$ ranges over infinitely many neighbours of $v"$,
therefore  $\gamma_{\sigma'(n)}^{-1} \partial (\G_v)$ converges to the unique
point of $\partial (\G_{v"})$ which we call $\zeta'$.   Therefore, the
convergence is locally uniform away from $\zeta'$, what we wanted to prove.
\endproof

As an immediate corollary of the two previous lemmas, we have:

\begin{coro}\label{lem;convgpe}
 With the previous notations, the group $\Gamma$ is a convergence group on $M$
(cf Definition \ref{def;convergence}).
\end{coro}

\begin{lemma}\label{lem;conic1}
  Every point in $\partial T \subset M$ is a conical limit point for $\Gamma$ in
$M$.
\end{lemma}

 Let $\eta \in \partial T$. Let $v_0$ a vertex in $T$ with a sequence
$(\gamma_n)_{n\in \mathbb{N}}$ of elements of $\Gamma$ such that $\gamma_n v_0$
lies on the ray $[v_0, \eta)$, converging to $\eta$.

 By Lemma \ref{lem;gdetrans}, after possible extraction of subsequence, there
is a point $\xi^{+} \in M$, and for all $\xi \in M$, except possibly one in
$\partial T$,  we have $\gamma_n^{-1} \zeta \to \xi^{+}$. Note that, in
particular, we have $\gamma_n^{-1} \partial(\G_{v_0}) \to \xi^{+}$. By
multiplying each $\gamma_n$ on the right by elements of $\G_{v_0}$, we can
assume that $\xi^{+}$ is not in $\partial(\G_{v_0})$, and we still have
$\gamma_n v_0$ lying on the ray $[v_0, \eta)$, converging to $\eta$.

Now it is enough to show that $\gamma_n^{-1}.\eta$ does not converge to
$\xi^{+}$. But $v_0$ is always in the ray $[\gamma_n^{-1} v_0,
\gamma_n^{-1}\eta)$. Therefore, if $\gamma_n^{-1}\eta \to \xi^{+}$, this implies
that $\xi^{+}$ is in $\partial(\Gamma_{v_0})$, which is contrary to our choice
of $(\gamma_n)_{n \in \mathbb{N}}$. \endproof

\begin{lemma}\label{lem;conic2}
  Every point in $\Omega/_{\sim}$ which is image by $\pi$ of a conical limit
point in a vertex stabilizer's boundary, is a conical limit point for $\Gamma$.
\end{lemma}

Such a point is in $\partial (\G_v)$ for some vertex $v$, and it is a conical limit point in
$\partial (\G_v)$  for $\G_v$.
Therefore it is a conical limit point in $M$ for $\G_v$ (see the remark
(i) in section 1), hence for $\G$. \endproof

\begin{lemma}\label{lem;parab}
  Every  point in $\Omega/_{\sim}$ which is image by $\pi$ of a bounded
parabolic point in a vertex stabilizer's boundary, is a bounded
parabolic point for $\Gamma$. The maximal parabolic subgroup associated is the
image in $\Gamma$ of a parabolic subgroup of a vertex group.

\end{lemma}

Let $\xi$ be the image by $\pi$ of a bounded parabolic point in a vertex
stabilizer's boundary, let $D(\xi)$ be its domain, and $v_1,\dots,v_n$ the
(finite, by Proposition \ref{coro;acyl}) list of vertices in $D(\xi)$ modulo
the action of $Stab(D(\xi))$, with stabilizers $\G_{v_i}$. Let $P$ be the
stabilizer of $\xi$. It stabilizes also $D(\xi)$, which is a bounded subtree of
$T$. By the Serre fixed-point theorem, it fixes a point in $D(\xi)$, which can be chosen
to be a vertex, since the action is without inversion. Therefore, $P$ is a
maximal parabolic subgroup of a vertex stabilizer, and the second assertion of
the lemma is true. For each $i\leq n$ the corresponding maximal parabolic
subgroup $P_i$ of $\G_{v_i}$ is a subgroup of $P$, because it fixes $\xi$. But
for each $i\leq n$, $P_i$ is bounded parabolic in $\G_{v_i}$, and  acts properly
discontinuously co-compactly on $\partial(\G_{v_i})\setminus
\{\xi\}$.

 For each index $i\leq n$, we choose $K_i \subset \partial(\G_{v_i})\setminus
\{\xi\} $, a compact fundamental domain of this action. We consider also
$\mathcal{E}_i$ the set of edges starting at $v_i$ whose boundary intersects
$K_i$ and does not contain $\xi$. Let $e$ be an edge with only one vertex in
$D(\xi)$, and $v_i$ be this vertex. As $K_i$ is a fundamental domain for the
action of $P_i$ on $\partial(\G_{v_i})\setminus
\{\xi\}$, there exists $p\in P_i$ such that $\partial (\G_e) \cap p K_i \neq
\emptyset$. Therefore, the set of edges $\bigcup_{i\leq n} P \mathcal{E}_i $
contains every edge with one and only one vertex in $D(\xi)$.

For each $i\leq n$, let $\mathcal{V}_i$ be the set of vertices $w$ of the tree
$T$ such that the first edge of
$[v_i,w]$ is in $\mathcal{E}_i$, and let $\overline{ \mathcal{V}_i }$ be its
closure in $T\cup \partial T$.  Let $K'_i$ be the subset of $M$
consisting of the points whose domain is included in $\overline{ \mathcal{V}_i }$.
  As a sequence of
 points in the boundaries of the stabilizers of distinct edges in
$\mathcal{E}_i$ has only accumulation points in $K_i$, the set  $K''_i=K_i\cup
K'_i$ is compact. Hence $\bigcup_{i\leq n} K''_i$ is a compact space not
containing $\xi$, and because $\bigcup_{i\leq n} P \mathcal{E}_i $
contains every edge with one and only one vertex in $D(\xi)$, the
union of the translates of $\bigcup_{i\leq n} K''_i$ by $P$  is $M \setminus
\xi$. Therefore, $P$ acts properly discontinuously co-compactly on $M\setminus
\xi$. \endproof

We can summarize the results of this section:

\begin{theo}[Dynamic of $\Gamma$ on $M$]\label{theo;dynamic}\qua

Under the conditions of Theorem 0.1, and with the previous notations, the group
$\Gamma$ is a geometrically finite convergence group on $M$.

The bounded parabolic points are the images by $\pi$ of bounded parabolic
points, and their stabilizers are the images, and their conjugates, of maximal
parabolic groups in vertex groups. \end{theo}

We are now able to prove our main theorem.

\proof[Proof of Theorem 0.1] The two first cases are  direct consequences
of Theorem \ref{theo;dynamic} and of Theorem \ref{theo;asli}. The maximal
parabolic subgroups are given by  Lemma \ref{lem;parab}.

Cases (3) and $(3')$ can be deduced as follows. Let $\Gamma = G_1 *_P G_2$,
where $P$ is \emph{maximal} parabolic in $G_1$ and parabolic in $G_2$. If
$\widetilde{P}$ is the maximal parabolic subgroup of $G_2$ containing $P$, one
has  $\Gamma = (G_1 *_P \widetilde{P})*_{\widetilde{P}} G_2$. One can apply
successively the second and the first case of the theorem to get the relative
hyperbolicity of $\Gamma$. For the last case, If $\Gamma = G*_P$, then one can
write $\Gamma = (G*_P P')*_{P'}$, where $P'$ is as in the statement, and apply 
consecutively the second and first case of the theorem. The acylindricity or the
last HNN-extension is given by the fact that the images of $P'$ in the group
$(G*_P P')$, are maximal parabolic subgroups not in the same conjugacy class.
\endproof

\section{Relatively Hyperbolic Groups and Limit Groups}

In our combination theorem, the construction of the boundary helps us to get
more information. For instance,  we get an independent proof, and an extension to the
relative case, of a theorem of I~Kapovich \cite{Kapo} for hyperbolic groups.

\begin{coro}
  If $\Gamma$ is in Case (1) of Theorem 0.1, the vertex groups embed as
fully quasi-convex subgroups in $\Gamma$. \end{coro}

The limit set of the stabilizer of a vertex $v$ is indeed $\partial (\G_v)$.
As domains are finite (Proposition \ref{coro;acyl} and its remark), a point in
$M$ belongs to finitely many translates of $\partial (\G_v)$. \endproof

Finally, we study limit groups, introduced by Sela in \cite{SelaLimit}, in his
solution of the Tarski problem, as a way to understand the structure of the
solutions of an equation in a free group. We give the definition of limit groups
; it involves a Gromov-Hausdorff limit. Here, we do not discuss the existence of
such a limit, but we advise the reader to refer to Sela's original paper.

\begin{defi}(Limit groups, \cite{SelaLimit})

Let $G$ be a finitely generated group, with a finite generating family $S$, and
$\gamma= (\gamma_1 \dots \gamma_k)$ a prescribed set of $k$ elements in $G$. Let
$F$ be a free group of rank $k$ with a fixed basis $a= (a_1 \dots a_k)$, and let
$X$ be its associated Cayley graph (it is a tree). Let $H(G, F ; \gamma,a)$ be
the set of all the homomorphisms of $G$ in $F$ sending $\gamma_i$ on $a_i$. Each
element of  $H(G, F ; \gamma, a)$ naturally defines an action of $G$ on $X$. Let
$(h_n)_{n\in \mathbb{N}}$ be a sequence of homomorphisms in distinct conjugacy
classes, and let us rescale $X$ by a constant $\mu_n =\min_{f\in F}  \max_{g\in
S} (d_X (id,f h_n(g) f^{-1} ))$ to get the pointed tree $(X_n,x_n)$, whose
base point $x_n$ is the image of a base point in $X$. There is a subsequence
such that $(X_{\sigma (n)},x_{\sigma (n)})$ converges in the sense of
Gromov-Hausdorff, and let $(X_{\infty},y)$ be the real tree that is the
Gromov-Hausdorff limit, on which the group $G$ acts. Let $K_{\infty}$ be the
kernel of this action (the elements of $G$ fixing every point in $X_{\infty}$).
We say that the quotient $L_{\infty} = G/K_{\infty}$ is a limit group.
\end{defi}

An important property of limit groups is an accessibility
theorem, proven by Sela. Every limit group has a height: limit groups of height $0$ are the finitely generated
torsion-free abelian groups, and every
limit groups of height $n>0$ can be constructed by
finitely many free products, HNN-extensions or amalgamations of limit groups of
height at most $(n-1)$, over cyclic groups (this is a consequence of Theorem 4.1
in \cite{SelaLimit}). We need only this fact, and the fact that every abelian
subgroup of a limit group is contained in a \emph{unique} maximal abelian
subgroup (Lemma 1.4 in \cite{SelaLimit}). Limit groups are known to enjoy many
more powerful properties, therefore, one can hope that a similar argument than
ours would work for a wider class of groups. We now establish
properties of acylindricity, which can be also found in \cite{SelaLimit}.

\begin{lemma}

One can choose the accessibility splitings of a limit group to be acylindrical.
Moreover, the edge group of a spliting involved is a maximal abelian subgroup,
and malnormal, in at least one of the adjacent vertex groups.

\end{lemma}

 If an amalgamation $A*_Z B$ is involved in the accessibility, then
 the subgroup $Z$ is maximal abelian in either $A$ or $B$, since it is a
property of limit groups that it has to belong to a \emph{unique} maximal
abelian subgroup of the ambiant group. In particular $Z$ is malnormal in either
$A$ or $B$, since if it was not, a proper subgroup would be in two distinct
 maximal abelian subgroups. Hence the amalgamation is 3--acylindrical.

 If an HNN-extension $A*_Z$ is involved, then, let $Z_1$ and $Z_2$ be the two
 images of $Z$ in $A$ in the extension, and let $t$ be a generator of the
loop of the graph of group.
 If $Z_1$ is not maximal abelian, let  $a_1$ be an
element, not in $Z_1$, and in the unique maximal abelian  subgroup of $A$
containing $Z_1$. If one conjugates $a_1$ by $t$, one gets an element of
$A *_Z$ not in $A$, that commutes with $Z_2$.
Then $Z_2$ is maximal abelian in $A$, as if it was not, it
would not be in a unique maximal abelian subgroup of $A*_Z$. Therefore, as in
the case of amalgamations, we see that either $Z_1$ or $Z_2$ has to be maximal
abelian, and therefore malnormal.
 Now, unless $Z_1 = Z_2$, we see that they cannot intersect non-trivially,
 because they would span a larger maximal abelian subgroup, contradicting what
we just proved. Therefore, the  HNN-extension is 2--acylindrical. Finally, if
$Z_1 = Z_2$, note that $A*_Z = (A*_Z Z)*_Z =A*_Z (Z*_Z) =  A *_Z Z^2$, which is
a previous case. \endproof

From this accessibility, Sela deduces that limit groups are exactly the finitely
generated $\omega$--residually free groups: these are the groups such that, for
every finite family of non-trivial elements,  there exists a morphism in a free
group that is non trivial on each of these elements.

 We will need the general fact:
 \begin{lemma}\label{lem;one_more_parab}

  Let $(G,\mathcal{G})$ be a relatively hyperbolic group, and let
 $Z$ be a non parabolic infinite cyclic subgroup of $G$ which is its
 own normalizer. Let $\mathcal{Z}$ be the set of conjugates of $G$. Then
 $(G,(\mathcal{G}\cup \mathcal{Z}))$ is a relatively hyperbolic group.
 \end{lemma}

To see this, note that the space $M$ obtained from $\partial (G)$ by
identifying for each conjugate of
$Z$, the two points of its limit set to a point, is
Hausdorff because the sequence of the diameters of the preimages in $\partial
(G)$ of any sequence of points in $M$ tends to zero (this is a
consequence of Proposition \ref{prop;diam}, for instance).
 Therefore, $M$ is a compact
metrisable space, on which the group $G$ acts as a convergence group. The images
in $M$ of bounded parabolic points of $\partial (G)$ are still bounded parabolic
points, with same stabilizers.
 If $\xi\in M $ is the image of a conical limit point, not in the
limit set of some conjugate of $Z$, there is a sequence $(g_n)$ in $G$, and $a$
and $b$ distinct points of $\partial (G)$  such that $g_n \xi \to a$ and $g_n
\zeta \to b$ for all other $\zeta$. If $a$ and $b$ map to the same point in $M$,
then they are in the limit set of a same conjugate $Z'$ of $Z$.
We assumed that $\xi$ is not in the limit set of $Z'$. Then by multiplying
the $g_n$ by sufficiently large elements $z_n$ of $Z'$ we  would get that
$z_n g_n \xi \to a$, $z_n g_n b \to b$ and for a sequence of points $a_n$
tending to $a$ more slowly than $g_n \xi$, $z_n g_n a_n \to c$ a point in a
fundamental domain of $Z'$ acting on $\partial(B)\setminus \Lambda Z'$. In
particular, this violates the convergence property. Therefore the images in $M$
of $a$ and $b$ are distinct.  Hence, the sequence
$g_n$ and the images of $a$ and $b$ in $M$ show that $\xi$ is a conical limit
point.

  If $\xi\in M $ is the image of the limit set of $Z$
(which consists of two loxodromic fixed points), then its stabilizer is the
normalizer of $Z$, that is $Z$ itself. As the cyclic group $Z$ acts
co-compactly on the complement of its limit set in $\partial (G)$, (this is a
consequence of the fact that $Z$ acts as a convergence group on $\partial G$
fixing the two points in its limit set), we see that $\xi$ is a bounded
parabolic point in $M$. Similar fact is true for every  conjugate of $Z$. All
this together proves the relative hyperbolicity, by Theorem \ref{theo;asli}.
\endproof

\begin{theo}{\rm [Theorem 0.3]}

 Every  limit group is hyperbolic relative to the family of
its maximal non-cyclic abelian subgroups. 
\end{theo}

We argue by induction on the height. It is obvious for groups of height $0$.
Consider an  HNN extension $A*_Z$ or an amalgamation $A*_Z B$, with $A$ and $B$
of height at most $(n-1)$, $Z$ cyclic. If $Z$ is trivial or has cyclic
centralizer in the amalgamation, both of its images in the vertex group(s)  are
fully quasi-convex, because it has trivial intersection with
every non-cyclic abelian subgroup. Hence, the first case of the combination
theorem 
 gives the result.

Assume now that $A$ contains a maximal
non-cyclic abelian subgroup containing $Z$.
We consider the case of an amalgamation $A*_Z B$, the case of an HNN-extension
being similar. Let $\{P_i\}$ be the set of maximal parabolic subgroups  of $B$;
each $P_i$ is a non-cyclic abelian group. From the discussion on the
accessibility, we know that the group $Z$ is a maximal cyclic subgroup of $B$
not intersecting any of the $P_i$, and is malnormal in $B$.  In particular, it
is fully quasi-convex in $B$, and we note ${Z}_i$ the set of conjugates of $Z$.
 From Lemma \ref{lem;one_more_parab}, we have that $B$ is hyperbolic
relative to $\{P_i\} \cup \{Z_i\}$.

We can apply the third case of Theorem \ref{theo;comb},
this gives that $A*_Z B$
 is hyperbolic relative to its
maximal non-cyclic abelian subgroups,
and this ends the proof for amalgamations.

The proof is similar in the case of an HNN-extension, using the case $(3')$ of the
combination theorem, instead of the third case. \endproof

The next proposition was suggested by G~Swarup (see also \cite{Sw}). It was
already known that every finitely generated subgroup of a limit group is itself a limit group
(it is obvious if one thinks of $\omega$--residually free groups).

\begin{prop}[Local quasi-convexity]\label{prop;coherence}\Gap

  Every finitely generated subgroup of a limit group is quasi-convex (in the sense of Definition \ref{def;full-q-c}).

  \end{prop}

Again, we argue by induction on the height of limit groups.

The result is  classical for free groups, surface groups, and abelian groups.
Assume now that the property is true for $A$ and $B$, and consider $\Gamma =
A*_Z B$, and $H$ a finitely generated subgroup of $\Gamma$. $H$ acts on the
Serre tree $T$ of the amalgamation. In particular it acts on its minimal invariant subtree.
As  a consequence of the fact that $H$ is finitely generated, the quotient of
this tree is finite. Moreover, as the edge groups are all cyclic or trivial, $H$
intersects each stabilizer of vertex along a finitely generated subgroup.
Therefore, one gets a spliting of $H$ as a finite graph of groups, the vertex
groups of which are finitely generated subgroups of the conjugates of $A$ and
$B$, and with cyclic or trivial edge groups. As they are finitely generated, and
by the induction assumption, the vertex groups are quasi-convex in the
conjugates of $A$ and $B$, and their boundaries equivariantly embed in the
translates of $\partial A$ and $\partial B$. We can apply our combination
theorem on this acylindrical graph of groups, and as the Serre tree of the
splitting of $H$ embeds in the Serre tree of the splitting of $\Gamma$, its
boundary equivariantly embeds in $\partial T$. Thus, $H$ is a geometrically
finite group on its limit set in the boundary of $\Gamma$, hence it is
quasi-convex in $\Gamma$. \endproof

The Theorem \ref{theo;Howson} (Howson property for limit groups)  was motivated by a discussion with
G~Swarup. To prove it, we
first prove the Proposition \ref{prop;quasiconvabeliens}, inspired by some results
in \cite{Swaruplimitsets}:
we study the intersection of (not necessarly fully) quasi-convex subgroups.

This study completes the work of I~Kapovich, who proved the Howson property for
limit groups without any non-cyclic abelian subgroup (see \cite{Kapofaux} and
\cite {Kapoerratum}).

\begin{theo}\label{theo;Howson}
   Limit groups have the Howson property: the intersection of two finitely
generated subgroups is finitely generated.
\end{theo}

We postpone the proof, because we need the following:

\begin{prop}[Intersection of quasi-convex
subgroup]\label{prop;quasiconvabeliens}\Gap

Let $\Gamma$ be a relatively hyperbolic group, with only  \emph{abelian}
parabolic subgroups. Let $Q_1$ and $Q_2$ be  two
quasi-convex subgroups. Then $Q_1 \cap Q_2$ is quasi-convex. Moreover,
$\Lambda(Q_1\cap Q_2)$ differs from $\Lambda (Q_1) \cap \Lambda (Q_2)$ only by
isolated points.
 \end{prop}

Let $Q_1$ and $Q_2$ be two quasi-convex  subgroups of $\Gamma$ and $Q = Q_1 \cap
Q_2$. The limit sets satisfy $\Lambda(Q) \subset \Lambda(Q_1) \cap
\Lambda(Q_2)$, and the action of $Q$ on $\Lambda(Q)$ is of convergence. As in
Proposition \ref{lem;interfullquasiconv}, the conical limit points in
$\Lambda(Q)$ are exactly the conical limit points in $\Lambda(Q_1)$ and in
$\Lambda(Q_2)$. We want to prove that each remaining point in $\Lambda(Q)$ is a
bounded parabolic point. Those points are among the parabolic points in both
$\Lambda(Q_1)$ and $\Lambda(Q_2)$, but it may happen that a parabolic point for
$Q_1$ and $Q_2$ is not in $\Lambda(Q)$.

However, it is enough to prove that, for all $p$, parabolic point for $Q_1$ and
$Q_2$, then the quotient $Stab_Q(p) \backslash (\Lambda(Q_1) \cap \Lambda(Q_2)
\setminus \{p\})$ is compact. Indeed, if we manage to do so, we would have
proven that $\Lambda(Q)$ differs from $\Lambda(Q_1) \cap \Lambda(Q_2)$ only by
isolated points: the parabolic points for $Q_1$ and $Q_2$ whose
stabilizer in $Q$ is finite. Such a point $p$ is isolated, because the statement
above implies that $(\Lambda(Q_1) \cap \Lambda(Q_2) \setminus \{p\})$ is
compact. Therefore, Proposition \ref{prop;quasiconvabeliens}  follows from the
general lemma:

\begin{lemma}
Let $G$ be  a finitely generated   abelian group, acting properly
discontinuously on a space $E$. Assume that $G$ contains two subgroups, $A$ and
$B$, such that $G=AB$. If $A$ acts on $X\subset E$ with compact quotient, and if
$B$ acts similarly on $Y\subset E$, then $A\cap B$ acts properly discontinuously
on $X\cap Y$, with compact quotient. \end{lemma}

The  only thing that needs to be checked is that the quotient is compact. Let
$K_A \subset X$ be a compact fundamental domain for $A$ in $X$, and $K_B$
similarly for $B$ in $Y$. For all $a\in A$ such that $aK_A \cap Y \neq
\emptyset$, there exists $b \in B$ such that  $aK_A \cap bK_B \neq \emptyset$.
As $K_A$ and $K_B$ are compact, and since the action of $(A+B)$ is properly
discontinuous, there are finitely many possible values in $G$ for $a^{-1}b$, with $a$ and $d$
satisfying $aK_A \cap bK_B \neq \emptyset$. Therefore, for all such $a$ and $b$, there
exists a word $w$ written with an alphabet of generators of $G$ consisting of
generators of $A$ and generators of $B$, of length bounded by a number $N$ neither
depending on $a$ nor on $b$, such that, in $G$, $w=a^{-1}b$.
Using abelianity of
the group $G$, we can gather the letters in $w$ in order to get a new word of
same length, concatenation of two smaller ones: $w'=w_Aw_B$ with $w_A \in A$
and $w_B \in B$, and still, in $G$,  $w'=a^{-1}b$. Now we see that $aw_A =
b(w_B)^{-1}$, and therefore $aw_A \in (A\cap B)$. If we set $K= (\bigcup_{|w_A|\leq
N}w_AK_A )\cap Y$, which is compact, we have just shown that $(A\cap B)K$ covers
$X\cap Y$. That is that we have proven the lemma. \endproof

Now we can prove the Howson property.

\proof[Proof of Theorem \ref{theo;Howson}]
Two finitely generated subgroups of a limit group are quasiconvex by Proposition
\ref{prop;coherence}, therefore, by Proposition \ref{prop;quasiconvabeliens},
the intersection is also quasiconvex. In particular, by remark {(ii)}
in section 1, it is finitely generated. \endproof

We finally give an application of \cite{Dah}. Following Bestvina \cite{Be}, we
say that a $\mathcal{Z}$--structure (if it exists) on a group is a minimal
(in the sense of $\mathcal{Z}$--sets) aspherical equivariant, finite dimensional
(for the topological dimension) compactification of a universal
cover of a finite
classifying space for the group, $E\Gamma \cup \partial(E\Gamma)$, such that
the convergence of a sequence $(\gamma_n p)_n$ to a point of the boundary
$\partial(E\Gamma)$ does not depend on the choice of the point $p$ in $E\Gamma$
(see \cite{Be}, \cite{Dah}).

\begin{theo}[Topological compactification]\label{theo;Zstruc}\Gap

Any limit group admits a $\mathcal{Z}$--structure in the sense of \cite{Be}.
\end{theo}

The maximal parabolic subgroups are isomorphic to some $\mathbb{Z}^d$, and
therefore admits a  finite classifying space with a $\mathcal{Z}$--structure (the
sphere that comes from the $CAT(0)$ structure). As limit groups are torsion
free, (Lemma 1.3 in \cite{SelaLimit}), the main theorem of \cite{Dah} can be
applied to give the result. \endproof

 We emphasize that this topological boundary needs not to be the one constructed
above: if the group contains $\mathbb{Z}^d$, the topological boundary contains
a sphere of dimension $d-1$.


 \end{document}